\documentclass{llncs}

\usepackage[pdftex]{graphicx}
\usepackage{epstopdf}        
     
\usepackage{makeidx}  
\usepackage{amsmath}
\usepackage{amsfonts,amssymb}
\usepackage{clrscode3e}        
\usepackage{multirow}

\newcommand{\dt}{\displaystyle}

\newcommand{\mbf}[1]{\protect\text{\boldmath$#1$}}                 

\begin{document}
\frontmatter          
\pagestyle{headings}  
\addtocmark{The set cover problem} 
\title{On the representativeness of approximate solutions of discrete optimization problems with interval objective function}
\titlerunning{On the representativeness of approximate solutions...}  
%
\author{Alexander Prolubnikov}
\authorrunning{Alexander Prolubnikov} 
%
\tocauthor{Alexander Prolubnikov}
\institute{Omsk State University, Omsk, Russian Federation\\
\email{a.v.prolubnikov@mail.ru}
}
\maketitle              

\begin{abstract}

We consider discrete optimization problems with interval uncertatinty of objective function coefficients. The interval uncertainty models measurements errors. A pos\-sible optimal solution is a solution that is optimal for some possible values of the coefficients. Pro\-ba\-bi\-li\-ty of a possible solution is the probability to obtain such coefficients that the solution is optimal. Similarly we define the notion of a possible approximate solution with given accuracy and probability of the solution. A possible approximate solution is an approximate solution that is obtained for some possible values of the coefficients by some fixed approximate algorithm, e.g. by the greedy algorithm. Pro\-ba\-bi\-li\-ty of a such solution is the probability to obtain such coefficients that the algorithm produces the solution as its output. We consider optimal or approximate possible solution un\-re\-pre\-sen\-ta\-ti\-ve if its probability less than some boundary value. The mean approximate solution is a possible approximate solution for midpoints of the coefficients intervals. The solution may be treated as approximate solution for exact values of the coefficients. We show that the share of individual discrete optimization problems instances with unrepresentative mean approximate solution may be wide enough for rather small values of error and the boundary value. The same is true for any other possible approximate solution: all of them may be unrepresentative.

\end{abstract}

\section*{Introduction}

The ultimate goal of applied scientific investigations is a prediction of experimental results. This concerns such mathematical discipline as discrete optimization too. Discrete optimization's methods use measurement results as input data for its models. The models include {\it objective function} (or {\it cost function} or {\it loss function}) which sets such characteristic of the modeled phenomenon that it is needed to minimize (or maximize) its value. We may characterize our preferences using the notion of costs. The costs equals to losses that we shall have choosing some possible solution.

In the most general form, we may formulate discrete optimization problem that we consider (we shall call it as {\it DO-problem} further) in the following form.

\smallskip

\noindent{\bf DO-problem (I).} Let $E\!=\!\{e_1,\ldots,e_n\}$. $c(e)\!>\!0$ is a {\it cost} of $e\!\in\!E$, $c_i\!=\!c(e_i)$. A binary vector $x\!=\!(x_1,\ldots,x_n)$ defines the set $E_x\!\subset\!E$: $x_i\!=\!1$, if $e_i\!\in\! E_x$, and $x_i\!=\!0$, if $e_i\!\in\! E\setminus E_x$. The set $\mathcal{D}$ of feasible solutions is given. We need to find $\check{x}\!\in\!\mathcal{D}$ such that $$f(\check{x},c)=\min\limits_{x\in\mathcal{D}}f(x,c),$$
where
\begin{equation}\label{1}  
f(x,c)=\sum\limits_{e\in E_x}c(e)=\sum\limits_{i=1}^nc_ix_i.
\end{equation}

\smallskip

\noindent A lot of DO-problems on graphs and hypergraphs may be formulated as problems of the form (I). Namely, the set $E$ may be considered as the set of edges of the graph $G\!=\!\langle V, E\rangle$, while the set $\mathcal{D}$ may be considered as a set of some subgraphs of $G$. We may consider the set of binary vectors associated with subgraphs of specified form as the set of feasible solutions $\mathcal{D}$. E.g., it may be the set of paths which connects some two graph vertices or the set of its spanning trees or the set of its Hamiltonian cycles, matchings in graphs, cuts and so on. Not only problems on graphs (hypergraphs) may be formulated as DO-problems of the form (I). For example, the knapsack problem may be formulated this way too. 

Various applied problems may be stated in a such way. But since, generally, exact measurements of DO-problems' parameters  quite often is impossible, we have uncertain but not exact values of them. Also, absence of information or its variability with time may be the cause of the uncertainty. It is a frequent situation that the only reliable information that we have on parameters is intervals of its possible values. Additionally, we suppose the uniform distribution on the intervals as the most uninformative probability distribution.

A {\it scenario} is a vector of values from given intervals of possible coefficients of a cost function (\ref{1}). A {\it possible optimal solution} of a DO-problem with interval cost function (we shall call it {\it IDO-problems, i.e. interval DO-problem}) is a feasible solution which is optimal for some scenario. Similarly, a {\it possible approxi\-mate solution with accuracy $\alpha$} of an IDO-problem is a such $\tilde{x}\!\in\!\mathcal{D}$ that, for some scenario $c$, it holds that $f(\tilde{x},c)\le\alpha f(\check{x},c)$ for some fixed $\alpha\!>\!0$. Fixing the value of $\alpha$, we shall call it {\it possible approximate solution} for short.

We call the {\it mean approximate solution} the possible approximate solution that we obtain for scenario with components which are midpoints of coefficients' intervals. Let us denote it as $\tilde{x}_{\mu}$. 

Considering the interval to which multiple measurements of some parameter are belong, it is often the case when reseachers use the midpoint of it to estimate the value of the parameter and use it as its actual value for computations. As a result of implementation of a such approach we obtain the only one possible solution of the problem with uncertain parameter. For example, having DO-problem with uncertain cost function coefficients, we may obtain only the mean optimal solution or the mean approximate solution and we may treat it as a valid solution for~the~problem. Further we show that this approach is not justified quite frequently. 

Suppose we use some approximate algorithm to solve some DO-problem, e.g. the greedy algorithm.
{\it Probability of a possible approximate solution} $\tilde{x}$ is equal to probability of~obtaining such a scenario that we obtain $\tilde{x}$ as a result of the algorithm operating for the scenario. For the case when probability of a possible (optimal or approximate) solution less than some given boundary value, we call the solution as {\it un\-re\-presen\-tative}.  We use probability of a possible solution as~a~measure of its representativeness. The solution is {\it representative} otherwise. We call this boundary value as {\it boundary of representativeness} and denote it as $b$.

Generally, considering NP-hard DO-problems, we may have such instances of the the problems that we cannot find optimal solution for reasonable time even though these are instances with not more than hundreds variables. So we use approximate algorithms with guaranteed accuracy which operates in polynomial time. Greedy algorithms are examples of such algorithms. 

When an instance of a DO-problem turns into an instance of an IDO-problem under assumption of parameters uncertainty, we must answer to the following questions about the situation we have. 
\begin{itemize}
\item[a).] How many possible (optimal or approximate) solutions this instance have?
\item[b).] How much values of costs differ for these solutions?
\end{itemize}
There were used different ways to answer these questions ever since DO-problems with exact parameters has been studied. The notion of stability of~a~solution have been used usually to address it. An optimal solution called {\it stable} if it remains to be optimal while the values of costs are varied within some predefined intervals. We present the concept of  representativeness of a solution that is close to~the~concept of stability by its meaning. 
The approach we propose gives answers to the questions above. Having a possible solution of the instance, e.g. the mean solution, we may define how likely is it that we shall obtain this solution for arbitrary values of uncertain parameters. 

A possible approximate solution may be stable for a scenario $c$, i.e. there exists such $\delta\!>\!0$ that there are no another possible approximate solutions in~$\delta$-neighbourhood of $c\!=\!(c_1,\ldots,c_n)$ in $\mathbb{R}^n$, but it may be unrepresentative for some boundary of representativeness $b$ and for some $\delta$. Further we consider such an example of a DO-problem with a stable but unrepresentative mean approximate solution. Having interval coefficients we may have such sets of possible approximate solution that no one of them, including the mean approximate solution, will not be representative. This means that there will not be any solution that has some predefined accuracy with probability more than $b$.

Thus we consider any possible (optimal or approximate) solution as an representative of the whole set of all possible (optimal or approximate) solutions. If~we model repeating situations using some probabilty distribution on intervals of~possible costs, then the solution is more representative if we obtain it more often using some exact or approximate algorithm. If we need to take a decision only at~once, then one of the two solutions is more representative than another if	probability of its obtaining is greater than probability of obtaining another solution, i.e. the possible solution (whether optimal or approximate) is more representative than another if such scenario that the solution be optimal or $\alpha$-approximate for some $\alpha$ is more probable.

Besides, using the notion of representativeness we attain some shortness of~formulations. So, for example, if some value of $b$ is given, then we may briefly characterize the solution as unrepresentative instead of writing that "probability of the solution is less than $b$". Also, instead of~using the words "the share of instances for which the probability of the mean solution less than some boundary value $b$"\ we shall write "the share of instances with~unrepresentative mean solutions."

We show that all of the possible approximate solutions, including the mean approximate solutions of DO-problems, may be un\-re\-pre\-sen\-ta\-tive even for relatively small values of error rate measurements of costs and~for~rather small boundaries of representativeness. The situation that there is no any representative possible approximate  solutions occurs more often as dimentionality of the problem grows. So, in these situations, the set of all possible solutions (optimal or approximate), probabilities of the solutions, intervals of~the~solutions costs, probability distribution on the whole set of possible solutions costs~--- these are the factors that we may take into account trying to predict and minimize costs if we need to take a solution in the situation of interval uncertainty.\\

\section{Discrete optimization problems with interval costs}

\subsection{Interval uncertainty representation} 

We use bold fonts to represent interval values: $$\mbf{a}\!=\![\underline{\mbf{a}},\overline{\mbf{a}}]=\{a\in\mathbb{R}\ |\ \underline{\mbf{a}}\!\le\! a\!\le\!\overline{\mbf{a}}\},$$ where $\underline{\mbf{a}}$ is the lower bound of the interval and $\overline{\mbf{a}}$ is its upper bound, $\underline{\mbf{a}}\!\le\!\overline{\mbf{a}}$. $\mathbb{IR}$~denotes the set of all such intervals on $\mathbb{R}$. $\mathbb{IR}_+\!=\!\{\mbf{a}\!\in\!\mathbb{IR}\ |\ \underline{\mbf{a}}\!>\!0\}$.

For $\mbf{a}, \mbf{b}\!\in\!\mathbb{IR}$, $$\mbf{a}+\mbf{b}=[\underline{\mbf{a}}+\underline{\mbf{b}},\overline{\mbf{a}}+\overline{\mbf{b}}].$$ For $\mbf{a}\!\in\!\mathbb{IR}$, $\alpha\!\in\!\mathbb{R}_+$, $$\alpha\mbf{a}=[\alpha\underline{\mbf{a}},\alpha\overline{\mbf{a}}].$$ 

An interval vector $\mbf{a}$ is defined: $$\mbf{a}=(\mbf{a}_1,\ldots,\mbf{a}_n)=([\underline{\mbf{a}}_1,\overline{\mbf{a}}_1],\dots,[\underline{\mbf{a}}_n,\overline{\mbf{a}}_n]).$$ $\mathbb{IR}^n$ denotes the set of $n$-dimensional interval vectors with components from $\mathbb{IR}$. We consider scenarios from $\mathbb{IR}_+^n$.

For DO-problems with interval costs, we replace the cost function of the form (\ref{1}) with the interval cost function of the form
\begin{equation}\label{2}  
\mbf{f}(x,\mbf{c})=\sum\limits_{e\in E_x}\mbf{c}(e)=\sum\limits_{i=1}^n\mbf{c}_ix_i,
\end{equation}
\noindent where $\mbf{c}_i\!=\mbf{c}(e_i)\!\in\!\mathbb{IR}_+$ are intervals of possible cost values, $e_i\!\in\!E$. We use in (\ref{2}) the operations of addition and multiplication that we have defined. 

$\mbf{f}(x,\mbf{c})$ is the interval of possible costs for solution $x\!\in\!\mathcal{D}$. Since (\ref{2}) contains only the first powers of every variable and it contains them only at once, it follows from the~main theorem of~interval arithmetic \cite{SharyBook} that $$\mbf{f}(x,\mbf{c})=\{f(x,c)\ |\ c\in\mbf{c}\}=[f(x,\underline{\mbf{c}}),f(x,\overline{\mbf{c}})].$$ Note that interval $\mbf{f}(x,\mbf{c})$ contains all of the possible costs of $x$ for all possible scenarios $c\!\in\!\mbf{c}$. But, in practical applications, we need to know the subinterval of $\mbf{f}(x,\mbf{c})$ that contains costs of $x$ for such scenarios that $x$ is optimal or $\alpha$-approximate since else we may (should) choose another possible (optimal or approximate) solution instead of $x$.\\

\section{The interval greedy algorithm\\ for the set cover problem with interval costs}

Greedy algorithms utilize a rather common approach to obtain approximate solutions of~DO-problems of the form (I). Performing iterations of a greedy algorithm, we form the set $Gr\!\subset\!E$. We obtain the solution $x$ by choosing elements $e\!\in\!E$ one by one and putting them into~$Gr$ considering their costs $c(e)$ and other parameters of the instance. We implement this using some selection function to minimize the overall cost of the solution that we obtain adding element $e\!\in\!E$ into $Gr$. The algorithm stops operating when we obtain such $Gr$ that $Gr\!=\!E_x$ for some $x\!\in\!\mathcal{D}$. Let us call the solution that is obtained by the greedy algorithm for some scenario as a {\it greedy solution}.

If the set $\mathcal{D}$ is a~matroid, then the greedy algorithm gives an optimal solution \cite{Papadimitriu}. It has been proven that greedy algorithms are asymp\-to\-ti\-cal\-ly best polynomial algorithms for some DO-problems (\!\!\cite{Feige,FriezeSzpankowski} et al.). 

We shall consider {\it the set cover problem} (\!{\it SCP}) as an example of a DO-problem (I). There are given a set $\mathcal{U}\!=\!\{1,\ldots,m\}$ and such a collection $S$ of its subsets $S\!=\!\{S_1,\ldots, S_n\}$, $S_i\!\subseteq\!\mathcal{U}$, that $$\bigcup\limits_{i=1}^n S_i\!=\!\mathcal{U}.$$ We call a collection of sets $S'\!=\!\{S_{i_1},\ldots, S_{i_k}\}$, $S_{i_j}\!\in\! S$, a {\it cover} of  $\mathcal{U}$ if $$\bigcup\limits_{j=1}^{k}S_{i_j}\!=\!\mathcal{U}.$$ There are given costs $c_i\!=\! c\,(S_i)$, $c_i\!>\!0$ for $S_i\!\in\! S$, i.e. the vector $c\!\in\!\mathbb{R}_+^n$ of costs is given. Cost $c\,(S')$ of a collection of sets $S'\!=\!\{S_{i_1},\ldots, S_{i_k}\}$ is equal to the sum of costs of its elements: $$c(S')=\sum\limits_{j=1}^k c(S_{i_j}).$$ We need to find an {\it optimal cover} that is the cover with minimum cost.

Formulating SCP as a DO-problem of the form (I), we associate  the set $E$ with the collection $S$ of subsets from $\mathcal{U}$. For the problem, a binary vector $x\!\in\!\mathcal{D}$ of dimensionality $n$ defines a cover of $\mathcal{U}$: if $x_i\!=\!1$, then $S_i$ belongs to the cover, $x_i\!=\!0$ otherwise. 

Let us denote the cost of the solution $x\in\mathcal{D}$ as $c(x)$. For the case of interval costs, we define the interval of possible costs $\mbf{c}(x)$ of the solution $x\!\in\!\mathcal{D}$, i.e. $\mbf{c}(x)$ contains all of the costs of $x$ for all of the possible scenarios for which $x$ is a greedy solution. Note that for non-interval costs, we have $$c(x)=\sum_{\{i|x_i=1\}} c(S_i).$$ While for the case of interval costs, generally we have $$\mbf{c}(x)\neq\!\sum_{\{i|x_i=1\}} \mbf{c}(S_i),$$ since the interval $\mbf{c}(x)$ often may be refined from the value $\sum_{\{i|x_i=1\}} \mbf{c}(S_i)$ and we may have strict inclusion $$\mbf{c}(x)\subset\!\sum_{\{i|x_i=1\}} \mbf{c}(S_i).$$ 

Let us denote an optimal cover by $Opt$ and let $Cvr$ denote the cover produced by some algorithm Alg. Let $$\rho(\mbox{Alg})=\frac{c(Cvr)}{c(Opt)}.$$ SCP with real-valued costs is NP-hard \cite{GareyJohnson}. In \cite{Dinur}, it have been shown that $$\rho(\mbox{Alg})\!>\!(1-o(1))\ln m$$ for any polynomial time approximate algorithm Alg for SCP whenever $\mbox{P}\!\neq\!\mbox{NP}$. The results concerned with computational complexity of approximation for~unweighted case of the problem (\!\!\cite{Feige,RazSafra} et al.), i.e. when $c_i\!=\!1$, $i\!=\!\overline{1,n}$, are hold true for weighted case of SCP too.

Computational complexity of the greedy algorithm for SCP is of order $O(m^2n)$. Let us denote the cover produced by the greedy algorithm GrAlg as $Gr$. The~following logarithmical estimation of $\rho(Gr)$ holds for the algorithm \cite{Chvatal}: $$\rho(\mbox{GrAlg})=\dt\frac{c(Gr)}{c(Opt)}\le H(m)\le \ln m+1,$$ where $H(m)\!=\!\sum_{k=1}^m1/k$. Further we shall consider greedy solutions as possible approximate solutions of SCP with interval costs.

\smallskip

For the interval set cover problem ({\it ISCP}), there are given interval costs $\mbf{c}_i\!=\!\mbf{c}\,(S_i)$, i.e. an interval vector of possible scenarios $\mbf{c}\!\in\!\mathbb{IR}^n$ is given. As an optimal solution of an IDO-problem, we consider the {\it united solution set} of the problem that is the set $\Xi$ that contains possible optimal solutions for all of the scenarios in $\mbf{c}$: $$\Xi=\{x\in\mathcal{D}\ |\ (\exists c\in\mbf{c})(f(x,c)=\min\limits_{y\in\mathcal{D}}f(y,c))\}.$$ By now there is no algorithms to obtain $\Xi$ unless we not consider some type of~exhaustive search on scenarios from $\mbf{c}$. Such an algorithm may be constructed only for some special cases of operating with discrete intervals.

The notion of the IDO-problem's united approximate solution set is introduced in \cite{Prolubnikov}. The IDO-problem's {\it united approximate solution set} $\tilde{\Xi}_{\alpha}$ is a set that contains possible $\alpha$-approximate solutions for all scenarios in $\mbf{c}$: $$\tilde{\Xi}_{\alpha}=\bigr \{x\in\mathcal{D}\ \bigr |\ (\exists c\!\in\!\mbf{c})\ \bigl (f(x,c)\le\alpha\min\limits_{y\in\mathcal{D}}f(y,c)\bigr )\bigl \}.$$ Considering the greedy algorithm for SCP, $\tilde{\Xi}_{\alpha}$ contains greedy solutions for all of the scenarios, $\alpha\!=\!H(m)$. With that in~mind, further we denote the united approximate solution set simply as $\tilde{\Xi}$.

We solve ISCP using the interval greedy algorithm \cite{Prolubnikov,Prolubnikov2}. The algorithm is a~ge\-ne\-ra\-li\-za\-tion of the greedy algorithm for the case of interval costs. It gives $\tilde{\Xi}$, and, performing the algorithm's iterations, we obtain exact values of intervals $\mbf{c}(\tilde{x})$ of possible costs for $\tilde{x}\!\in\!\tilde{\Xi}$. If probability distributions on intervals of costs are additionally given, using the interval greedy algorithm we obtain the probabilities $\mbox{\sffamily P}(\tilde{x})$ for $\tilde{x}\!\in\!\tilde{\Xi}$. 

The approach presented in \cite{Prolubnikov,Prolubnikov2} may be applied to other IDO-problems that we may obtain for DO-problems of the form (I). Its computational complexity depends on cardinality of the set $\tilde{\Xi}$. The~complexity is exponential for the worst case. The interval greedy algorithm is polynomial if cardinatlity of $\tilde{\Xi}$ for an IDO-problem of dimensionality of $m$ is bounded by polynomial of $m$. 

\smallskip

Let us consider the following instance of ISCP with $m\!=\!7$, $n\!=\!11$. The collection $\mathcal{S}$ of subsets of $\mathcal{U}\!=\!\{1, 2, 3, 4, 5, 6, 7\}$ contains the sets $S_1\!=\!\{3, 5\}$, $S_2\!=\!\{4, 6\}$, $S_3\!=\!\{1,3\}$, $S_4\!=\!\{2, 3, 4\}$, $S_5\!=\!\{1, 5, 6\}$, $S_6\!=\!\{4, 5, 6\}$, $S_7\!=\!\{1, 4, 6, 7\}$, $S_8\!=\!\{1, 3, 4, 6\}$, $S_9\!=\!\{2, 4, 5, 7\}$, $S_{10}\!=\!\{1, 3, 6, 7\}$, $S_{11}\!=\!\{1, 2, 4, 6\}$. The mean values of costs are given by vector $c_{\mu}\!=\!(119, 117, 124, 135, 128, 130, 143, 144, 144, 142, 141)$. Suppose that the relative error of~mea\-su\-re\-ments of costs is not greater than $5\%$. That is to say that $\mbf{c}_i\!=\![c_{\mu,i}-\delta,c_{\mu,i}+\delta]$, where $c_{\mu,i}$ is the $i$-th component $c_{\mu}$, $\delta\!=\!5$. 

Let us denote the possible approximate solutions as $\tilde{x}^{(i)}$, $i\!=\!\overline{1,k}$, $k\!=\!|\tilde{\Xi}|$. The set $\tilde{\Xi}$ for the problem is listed below. Also we give the intervals of possible values of $c(\tilde{x}^{(i)})$ and probabilities $\mbox{\sffamily P}(\tilde{x}^{(i)})$.
\begin{itemize}
\item[1)] $\tilde{x}^{(1)}\!=\!(1,0,0,1,0,0,1,0,0,0,0)$, $\mbf{c}(\tilde{x}^{(1)})\!=\![382,410]$, $\mbox{\sffamily P}(\tilde{x}^{(1)})\!=\!0.1542$;
\item[2)] $\tilde{x}^{(2)}\!=\!(1,0,0,0,0,0,1,0,1,0,0)$, $\mbf{c}(\tilde{x}^{(2)})\!=\![391,410]$, $\mbox{\sffamily P}(\tilde{x}^{(2)})\!=\!0.0007$;
\item[3)] $\tilde{x}^{(3)}\!=\!(1,0,0,0,0,0,1,0,0,0,1)$, $\mbf{c}(\tilde{x}^{(3)})\!=\![390,410]$, $\mbox{\sffamily P}(\tilde{x}^{(3)})\!=\!0.1342$;
\item[4)] $\tilde{x}^{(4)}\!=\!(0,0,0,0,0,0,0,1,1,0,0)$, $\mbf{c}(\tilde{x}^{(4)})\!=\![278,295]$, $\mbox{\sffamily P}(\tilde{x}^{(4)})\!=\!0.1172$;
\item[5)] $\tilde{x}^{(5)}\!=\!(0,0,0,0,0,0,0,0,1,1,0)$, $\mbf{c}(\tilde{x}^{(5)})\!=\![278,293]$, $\mbox{\sffamily P}(\tilde{x}^{(5)})\!=\!0.3166$;
\item[6)] $\tilde{x}^{(6)}\!=\!(1,0,0,0,0,0,0,0,1,0,1)$, $\mbf{c}(\tilde{x}^{(6)})\!=\![389,417]$, $\mbox{\sffamily P}(\tilde{x}^{(6)})\!=\!0.0826$;
\item[7)] $\tilde{x}^{(7)}\!=\!(1,0,0,0,0,0,0,0,0,1,1)$, $\mbf{c}(\tilde{x}^{(7)})\!=\![387,417]$, $\mbox{\sffamily P}(\tilde{x}^{(7)})\!=\!0.1946$.
\end{itemize}
\noindent We obtain the mean approximate solution $\tilde{x}_{\mu}\!=\!\tilde{x}^{(1)}\!=\!(1,0,0,1,0,0,1,0,0,0,0)$ for scenario $c_{\mu}$, $f(\tilde{x}_{\mu},c_{\mu})=397$. The possible approximate solution $\tilde{x}_{\mu}$ is only one of the seven possible approximate solutions. Its probability equals to $0.1542$ while the most probable approximate solution is $\tilde{x}^{(5)}$ with probability of $0.3166$. For $\delta\!=\!5$, the mean solution $\tilde{x}_{\mu}$ is unrepresentative if $b\!\ge\!0.2$.

Note that the approximate solution $\tilde{x}_{\mu}$ is stable since the solution $\tilde{x}_{\mu}$ is the only greedy solution in $\delta$-neighbourhood of $c_{\mu}\!\in\!\mathbb{R}^n$ for $\delta\!\le\!0.5$. But, as the intervals of possible $c_i$ values grow, there will be several  possible greedy solutions. And the value of their costs will differ.\\ 

\section{The share of instances\\ with unrepresentative mean approximate solution}

\subsection{The sample of ISCP instances used in experiments}

To show that the share of ISCP instances with unrepresentative mean approximate solution may be wide enough, we generate the samples of ISCP instances using the algorthm that we present below. Implementing the algorithm, we generate sets $S_i\!\subseteq\!\mathcal{U}$ with random elements. We generate these sets until we obtain such a collection $S$ that every element of $\mathcal{U}$ is covered at least $q$ times by the sets from $S$, i.e. not less than $q$ sets from $S$ contain the element. The results of the experiments considered further were obtained for $q\!=\!3$. 

Costs of sets from $\mathcal{S}$ are defined as follows: 
\begin{equation}\label{4}
c_i:=100+10\cdot p_i+\eta,
\end{equation} 
\noindent where $p_i\!=\!|S_i|$, $\eta$ is an integer random value that is uniformly distributed on interval $[-5,5]$. Thus, if $S_i\!\in\!\mathcal{S}$ models some real life object then its cost $c_i$ is defined mainly by $p_i$. The random value $\eta$ models possible differencies of the objects which are close to each other by its base characteristics.

Having an instance of SCP with costs $c\!=\!(c_1,\ldots,c_n)\!\in\!\mathbb{R}_+^n$, we model inaccuracy of measuring of costs $c_i$ under assumption that it belongs to the interval $[c_i-\delta,c_i+\delta]$, i.e. we suppose $c_{\mu}\!=\!c$. If $\delta\!=\!5$ then modeled relative error is not greater than $5\%$. 

Generating the instances of ISCP this way, we obtain the sample of the population for given $m$ and $\delta$. The poulation consists of the pairs of the form:
\begin{itemize}
\item[1)] the set $\mathcal{P}^{(i)}(m,\delta)$ that consists of $1000$ ISCP instances generated by the algorithm; 
\item[2)] the vector $d^{(i)}(m,\delta)$ which is characterize the distribution of the instances accordingly to the value of probability $\mbox{\sffamily P}(\tilde{x}_{\mu})$: $$d^{(i)}(m,\delta)=(d_1^{(i)}(m,\delta),\ldots,d^{(i)}_{10}(m,\delta)),$$
\end{itemize}
\noindent where $d^{(i)}_k(m,\delta)$ is the share (in percent) of instances from $\mathcal{P}^{(i)}(m,\delta)$ for which $$\mbox{\sffamily P}(\tilde{x}_{\mu})\!\in\!((k\!-\!1)/10,k/10].$$
Further we show the computed average values of $d^{(i)}(m,\delta)$ for samples of the form
\begin{equation}\label{4}  
\{(\mathcal{P}^{(i)}(m,\delta), d^{(i)}(m,\delta))\}_{i=1}^{100},
\end{equation} $m\!=\!5,10,15,20$, $\delta\!=\!\overline{1,5}$. The average values of $d^{(i)}(m,\delta)$, $i\!=\!\overline{1,100}$, and its standard deviations are presented in Tables 1--8.

\bigskip

\fbox{\parbox{0.95\textwidth}{

\begin{center}
{\bf Generating ISCP instances.}
\end{center}

\noindent{\bf Input:} $m$, $q$, $\delta$.

\begin{itemize}
\item[1.] $i:=0$.
\item[2.] Do the following iterations while all elements of $\mathcal{U}$ are covered\\ less than $q$ times.
	\begin{itemize}
	   \item[2.1.] $i:=i+1$.
		 \item[2.2.] Generate cardinality $p_i$ of the set $S_i$:\\ $p_i$ is a random value uniformly distributed on $\mathcal{U}\!=\!\{1,\ldots,m\}$.
	   \item[2.3.] Assuming uniform distribution on the set $\mathcal{U}$,\\ randomly select $p_i$ elements from $\mathcal{U}$ to obtain $S_i$.
     \item[2.4.] Generate cost $c_i$: $$c_i:=100+10\cdot p_i+\eta,$$ where $\eta$ is an integer random value distributed uniformly\\ on the interval $[-5,5]$.
		 \item[2.5.] Form the interval cost $\mbf{c}_i$:\\ $\mbf{c}_i:=[c_i-\delta, c_i+\delta]$.	
	\end{itemize}
\end{itemize}

\smallskip

\noindent{\bf Output:} An ISCP instance with given $\mathcal{U}$, $\mathcal{S}\!=\!\{S_1,\ldots,S_n\}$, vector $\mbf{c}\!\in\!\mathbb{IR}_+^n$.\\

}}\\

Cardinality of the set of all ISCP instances that we may obtain using the~presented algorithm is substantially exceeds cardinality of the set of all SCP instances that we may generate in the curse of its implementation which grows exponentially with respect to $m$ \cite{DevittJackson}. Even though we cannot say that the sample of $100$ elements $\mathcal{P}^{(i)}(m,\delta)$ is large enough for some fixed $m$ and $\delta$, small dimensionalities of considered instances give us statistically stable results as it justified by the values of standard deviation of $d^{(i)}(m,\delta)$.\\

\subsection{The computational experiments}

\subsubsection{The share of ISCP instances\\ with unrepresentative mean approximate solution.} 

To perform the computational experiments considered below, we search for the set $\tilde{\Xi}$ for every ISCP instance that we generate using interval greedy algorithm. The results of the computations are presented in the tables below. 

The average values of the components $d^{(i)}(m,\delta)$ for the sample (\ref{4}), $\delta\!=\!\overline{1,5}$, are presented in Tables 1, 3, 5, 7. The standard deviations from the values are presented in Tables 2, 4, 6, 8. So, for example, the value in Table 7 for $m\!=\!20$, $\delta\!=\!5$ and $k\!=\!10$ means that, in average, there are $36.63\%$ of instances from $\mathcal{P}^{(i)}(m,\delta)$, $i\!=\!\overline{1,100}$, for which $\mbox{\sffamily P}(\tilde{x}_{\mu})\!\in\!(0.9,1]$.

\begin{table}[h!]
\caption{Average values of $d^{(i)}(5,\delta)$, $i\!=\!\overline{1,100}$.}
\begin{center}
\begin{tabular}{|c||c|c|c|c|c|c|c|c|c|c|}
\hline
$\delta\backslash k$ & $1$ & $2$ & $3$ & $4$ & $5$ & $6$ & $7$ & $8$ & $9$ & $10$\\
\hline
1 & 0 & 0.05 & 0.18 & 0.53 & 1.59 & 5.04 & 0.10 & 1.49 & 10.13 & 80.89 \\
\hline
2 & 0 & 0.06 & 0.35 & 1.20 & 3.64 & 4.04 & 3.65 & 7.70 & 7.39 & 71.96 \\
\hline
3 & 0 & 0.13 & 0.65 & 1.98 & 5.27 & 4.57 & 7.70 & 7.13 & 5.32 & 67.26 \\
\hline
4 & 0 & 0.20 & 1.11 & 3.31 & 6.20 & 7.39 & 7.21 & 6.17 & 5.62 & 62.77 \\
\hline
5 & 0 & 0.35 & 1.71 & 4.23 & 7.58 & 9.29 & 7.00 & 5.09 & 4.09 & 60.65 \\
\hline
\end{tabular}
\end{center}
\end{table}

\begin{table}[h!]
\caption{The standard deviation of $d^{(i)}(5,\delta)$.}
\begin{center}
\begin{tabular}{|c||c|c|c|c|c|c|c|c|c|c|}
\hline
$\delta\backslash k$ & $1$ & $2$ & $3$ & $4$ & $5$ & $6$ & $7$ & $8$ & $9$ & $10$\\
\hline
1 & 0 & 0.07 & 0.17 & 0.23 & 0.63 & 1.17 & 0.10 & 0.45 & 2.12 & 3.94 \\
\hline
2 & 0.01 & 0.08 & 0.20 & 0.42 & 0.88 & 0.97 & 1.07 & 1.58 & 1.66 & 5.50 \\
\hline
3 & 0.01 & 0.12 & 0.27 & 0.72 & 0.97 & 1.15 & 1.53 & 1.09 & 1.10 & 5.55 \\
\hline
4 & 0.02 & 0.13 & 0.39 & 0.94 & 1.27 & 1.28 & 1.44 & 1.27 & 1.00 & 6.03 \\
\hline
5 & 0.01 & 0.18 & 0.47 & 1.13 & 1.57 & 1.51 & 1.10 & 0.95 & 0.83 & 5.92 \\
\hline
\end{tabular}
\end{center}
\end{table}

\begin{table}[h!]
\caption{Average values of $d^{(i)}(10,\delta)$, $i\!=\!\overline{1,100}$.}
\begin{center}
\begin{tabular}{|c||c|c|c|c|c|c|c|c|c|c|}
\hline
$\delta\backslash k$ & $1$ & $2$ & $3$ & $4$ & $5$ & $6$ & $7$ & $8$ & $9$ & $10$\\
\hline
1 & 0 & 0.05 & 0.23 & 0.48 & 1.72 & 6.08 & 0.13 & 1.82 & 11.87 & 77.61 \\
\hline
2 & 0 & 0.07 & 0.42 & 1.25 & 3.55 & 5.09 & 3.73 & 9.48 & 8.84 & 67.57 \\
\hline
3 & 0.01 & 0.15 & 0.81 & 2.31 & 5.47 & 5.67 & 9.90 & 8.82 & 6.92 & 59.94 \\
\hline
4 & 0.01 & 0.23 & 1.32 & 3.84 & 7.09 & 8.57 & 9.09 & 8.46 & 8.87 & 52.53 \\
\hline
5 & 0.02 & 0.40 & 2.10 & 4.95 & 9.17 & 11.36 & 9.21 & 7.84 & 8.12 & 46.81 \\
\hline
\end{tabular}
\end{center}
\end{table}

\begin{table}[h!]
\caption{The standard deviation of $d^{(i)}(10,\delta)$.}
\begin{center}
\begin{tabular}{|c||c|c|c|c|c|c|c|c|c|c|}
\hline
$\delta\backslash k$ & $1$ & $2$ & $3$ & $4$ & $5$ & $6$ & $7$ & $8$ & $9$ & $10$\\
\hline
1 & 0.01 & 0.07 & 0.15 & 0.21 & 0.45 & 0.88 & 0.12 & 0.42 & 1.11 & 1.46 \\
\hline
2 & 0 & 0.08 & 0.20 & 0.36 & 0.69 & 0.79 & 0.71 & 1.00 & 1.01 & 1.72 \\
\hline
3 & 0.03 & 0.15 & 0.31 & 0.54 & 0.80 & 0.77 & 0.98 & 0.82 & 0.94 & 1.65 \\
\hline
4 & 0.03 & 0.17 & 0.42 & 0.69 & 0.83 & 0.91 & 1.08 & 1.34 & 1.20 & 2.65 \\
\hline
5 & 0.04 & 0.20 & 0.42 & 0.79 & 0.88 & 1.18 & 1.17 & 1.07 & 1.53 & 3.54 \\
\hline
\end{tabular}
\end{center}
\end{table}

\begin{table}[h!]
\caption{Average values of $d^{(i)}(15,\delta)$, $i\!=\!\overline{1,100}$.}
\begin{center}
\begin{tabular}{|c||c|c|c|c|c|c|c|c|c|c|}
\hline
$\delta\backslash k$ & $1$ & $2$ & $3$ & $4$ & $5$ & $6$ & $7$ & $8$ & $9$ & $10$\\
\hline
1 & 0 & 0.03 & 0.22 & 0.39 & 1.39 & 5.85 & 0.71 & 1.70 & 11.54 & 78.16 \\
\hline
2 & 0 & 0.05 & 0.35 & 1.08 & 3.07 & 5.46 & 4.06 & 9.28 & 9.45 & 67.19 \\
\hline
3 & 0 & 0.10 & 0.60 & 1.88 & 4.83 & 6.10 & 10.31 & 9.42 & 8.47 & 58.29 \\
\hline
4 & 0.01 & 0.18 & 1.08 & 3.32 & 6.91 & 8.72 & 9.91 & 9.64 & 10.54 & 49.67 \\
\hline
5 & 0.01 & 0.35 & 1.84 & 4.53 & 8.73 & 12.26 & 10.18 & 9.11 & 10.96 & 41.99 \\
\hline
\end{tabular}
\end{center}
\end{table}

\begin{table}[h!]
\caption{The standard deviation of $d^{(i)}(15,\delta)$.}
\begin{center}
\begin{tabular}{|c||c|c|c|c|c|c|c|c|c|c|}
\hline
$\delta\backslash k$ & $1$ & $2$ & $3$ & $4$ & $5$ & $6$ & $7$ & $8$ & $9$ & $10$\\
\hline
1 & 0 & 0.05 & 0.14 & 0.22 & 0.47 & 0.92 & 0.22 & 0.46 & 1.39 & 2.55 \\
\hline
2 & 0.01 & 0.09 & 0.25 & 0.41 & 0.85 & 0.87 & 0.99 & 1.11 & 1.15 & 3.98 \\
\hline
3 & 0.01 & 0.10 & 0.33 & 0.79 & 1.14 & 0.86 & 1.18 & 1.36 & 0.87 & 4.16 \\
\hline
4 & 0.03 & 0.15 & 0.51 & 1.07 & 1.32 & 1.26 & 1.11 & 1.06 & 1.73 & 3.59 \\
\hline
5 & 0.03 & 0.26 & 0.65 & 1.20 & 1.68 & 1.00 & 1.20 & 1.19 & 1.70 & 4.31 \\
\hline
\end{tabular}
\end{center}
\end{table}

\begin{table}[h!]
\caption{Average values of $d^{(i)}(20,\delta)$, $i\!=\!\overline{1,100}$.}
\begin{center}
\begin{tabular}{|c||c|c|c|c|c|c|c|c|c|c|}
\hline
$\delta\backslash k$ & $1$ & $2$ & $3$ & $4$ & $5$ & $6$ & $7$ & $8$ & $9$ & $10$\\
\hline
1 & 0 & 0.01 & 0.10 & 0.29 & 1.01 & 5.35 & 1.42 & 1.89 & 9.92 & 80.00 \\
\hline
2 & 0 & 0.04 & 0.19 & 0.79 & 2.43 & 5.69 & 4.12 & 9.37 & 9.46 & 67.91 \\
\hline
3 & 0 & 0.06 & 0.46 & 1.67 & 4.64 & 7.04 & 10.71 & 10.02 & 10.55 & 54.84 \\
\hline
4 & 0 & 0.11 & 0.78 & 2.74 & 6.46 & 9.48 & 10.87 & 11.10 & 12.94 & 45.52 \\
\hline
5 & 0.01 & 0.20 & 1.32 & 3.98 & 8.62 & 12.92 & 11.89 & 11.46 & 12.97 & 36.63 \\
\hline
\end{tabular}
\end{center}
\end{table}

\begin{table}[h!]
\caption{The standard deviation of  $d^{(i)}(20,\delta)$.}
\begin{center}
\begin{tabular}{|c||c|c|c|c|c|c|c|c|c|c|}
\hline
$\delta\backslash k$ & $1$ & $2$ & $3$ & $4$ & $5$ & $6$ & $7$ & $8$ & $9$ & $10$\\
\hline
1 & 0.01 & 0.03 & 0.11 & 0.22 & 0.40 & 0.91 & 0.37 & 0.37 & 1.72 & 2.58 \\
\hline
2 & 0 & 0.06 & 0.15 & 0.33 & 0.83 & 0.76 & 0.71 & 1.31 & 1.08 & 3.50 \\
\hline
3 & 0 & 0.07 & 0.25 & 0.57 & 1.32 & 0.90 & 1.31 & 0.98 & 0.88 & 3.64 \\
\hline
4 & 0.01 & 0.12 & 0.42 & 0.95 & 1.65 & 1.52 & 1.29 & 1.15 & 1.13 & 5.66 \\
\hline
5 & 0.02 & 0.15 & 0.58 & 1.25 & 2.09 & 1.44 & 1.40 & 1.24 & 1.31 & 5.85 \\
\hline
\end{tabular}
\end{center}
\end{table}

In particular, the results we present show that, having the boundary of re\-pre\-sen\-ta\-ti\-ve\-ness $b\!=\!0.9$ and the relative error that is not greater than $5\%$, the mean approximate solution is unrepresentative for more than half of the generated ISCP instances. It also shows that the share of ISCP instances with unrepresentative mean approximate value grows as dimensionality of the instances grows. It holds for all values of $b\!=\!k/10$, $k\!=\!\overline{1,9}$.\\

\subsection{The representativeness of possible approximate solutions\\ of the ISCP instances of higher dimensionalities}

Let $\mbox{\sffamily P}_{\mu}$ denotes the average value of probabilities $\mbox{\sffamily P}(\tilde{x})$ over all $\tilde{x}\!\in\!\tilde{\Xi}$ for an ISCP instance, $\mbox{\sffamily P}_{\max}$ denotes the maximum of $\mbox{\sffamily P}(\tilde{x})$ over all $\tilde{x}\!\in\!\tilde{\Xi}$. Let $\mbox{\sffamily MP}_{\mu}$ denotes the average value of $\mbox{\sffamily P}_{\mu}$ over all generated ISCP instances for the fixed value of $m$ and let $\mbox{\sffamily MP}_{\max}$ denote the average value of $\mbox{\sffamily P}_{\max}$ for them.

In the Table 9, we present the values of $\mbox{\sffamily MP}_{\mu}$ and $\mbox{\sffamily MP}_{\max}$ that were computed from 1000 of ISCP instances generated for every $m$. $\mbf{c}_i\!=\![c_i-\delta_i,c_i+\delta_i]$ for the instances. The values of $\delta_i\!=\!0.05c_i$ are defined for every $\mbf{c}_i$ individually. The data presented in Table 9 shows that the value $\mbox{\sffamily MP}_{\mu}$ decreases by one order of 10 as $m$ increases by the same order. In average, even the value of $\mbox{\sffamily P}_{\max}$ is less than $0.25$ for ISCP instances  with $m\!>\!100$, i.e. the mean approximate solution and any other possible approximate solution is unrepresentative for $b\!\ge\!0.25$.

\begin{table}[h!]
\caption{The values of $\mbox{\sffamily MP}_{\mu}$ and $\mbox{\sffamily MP}_{\max}$.}
\begin{center}
\begin{tabular}{|c||c|c|c|}
\hline
$m$ & $\mbox{\sffamily MP}_{\mu}$ & $\mbox{\sffamily MP}_{\max}$\\
\hline
100 & 0.0554 & 0.3675\\
\hline
200 & 0.0214 & 0.2636\\
\hline
300 & 0.0128 & 0.2235\\
\hline
400 & 0.0092 & 0.1984\\
\hline
500 & 0.0071 & 0.1794\\
\hline
600 & 0.0055 & 0.1656\\
\hline
700 & 0.0049 & 0.1603\\
\hline
800 & 0.0040 & 0.1498\\
\hline
900 & 0.0036 & 0.1420\\
\hline
1000 & 0.0033 & 0.1378\\
\hline
\end{tabular}
\end{center}
\end{table}

Note that there is a rather wide share of the instances with representative mean approximate solution for $b\!=\!0.9$ and for $m\!\le\!20$ considering Tables 1--8. For example, it is equal to $36.63\%$ for $m\!=\!20$ and $\delta\!=\!5$. But having the values of $m$ an order of magnitude larger, there is almost no such instances that the mean approximate or any other possible approximate solution is representative for $b\!=\!0.9$.\\

\subsection{The distribution of weights of possible approximate solution over all scenarios} 

Let us turn back to the ISCP instance that we considered earlier. For the instance, there are two possible optimal solutions $\check{x}^{(1)}\!=\!(0,0,0,0,0,0,0,1,1,0,0)$ and $\check{x}^{(2)}\!=\!(0,0,0,0,0,0,0,0,1,1,0)$. $f(\check{x}^{(i)},c)\!\in\![276, 296]$ for all $c\!\in\!\mbf{c}$, $i\!=\!\overline{1,2}$. $\check{x}^{(2)}$ is the mean optimal solution ($\check{x}_{\mu}\!=\!\check{x}^{(2)}\!=\!\tilde{x}^{(4)}$). Thus $\mbox{\sffamily P}(\check{x}_{\mu})\!=\!0.667$, i.e.~$\check{x}_{\mu}$ has maximum probability. Its cost equal to the mean value of all possible optimal solutions costs for all scenarios: $c(\check{x}_{\mu})\!=f(\check{x}_{\mu},c_{\mu})\!=\!286$. Possible deviation from~the~mean value is equal to $10$:  
\begin{equation}\label{5}
\max\limits_{{c\in\mbf{c} \atop \check{x}\in\Xi}}|f(\check{x}_{\mu},c_{\mu})-f(\check{x},c)|=10.
\end{equation} 

\begin{figure}[htbp!]
	\begin{center}
		\includegraphics[width=320px]{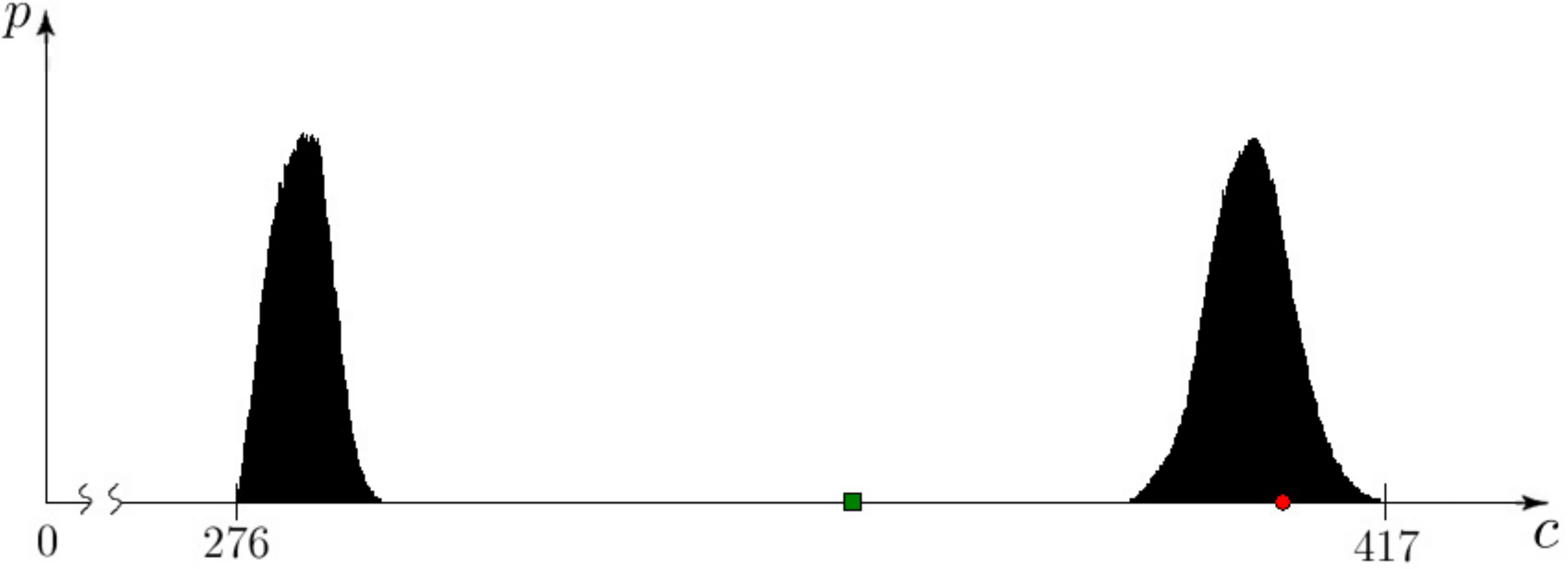}
       \caption {Hystogram of frequencies of costs for possible approximate solutions}
       \label{fig1}
   \end{center}
	\end{figure}
	
Hystogram of frequencies (axis $p$) of costs of possible approximate solutions (axis $c$) of the earlier consisdered ISCP instance is presented at Fig. 1. We build it based upon $10^6$ possible random scenarios that we generated using uniform distribution over intervals of possible costs for the instance. The hystogram approximates graph of costs distribution density. The cost of the mean approximate solution is equal to $350.51$. It is computed as the expected value as following: $$\sum\limits_{\tilde{x}\in\tilde{\Xi}}\mbox{\sffamily P}(\tilde{x})\cdot\frac{\underline{\mbf{c}}(\tilde{x})+\overline{\mbf{c}}(\tilde{x})}{2}\approx 350.51.$$ The value marked as a square on the axis $f$. The cost $f(\tilde{x}_{\mu},c_{\mu})$ of the mean approximate solution $\tilde{x}_{\mu}$ for midpoints scenario marked as a circle on the axis $f$, $f(\tilde{x}_{\mu},c_{\mu})\!=\!397$. Note that the set $\{f(\tilde{x},c)\ |\ (\tilde{x}\in\tilde{\Xi})(c\in\mbf{c})\}$ (i.e. the set of all of possible costs) is not connected for the instance.

We see a significant deviation of the mean approximate solution cost $f(\tilde{x}_{\mu},c_{\mu})\!=\!397$ from the expected cost of a possible approximate solution over all scenarios from $\mbf{c}$ and 
\begin{equation}\label{6}
\max\limits_{c\in\mbf{c} \atop \tilde{x}\in\tilde{\Xi}}|f(\tilde{x}_{\mu},c_{\mu})-f(\tilde{x},c)|=121.
\end{equation} 

Comparing (\ref{5}) and (\ref{6}), we may see that costs of possible approximate solutions may differ from the cost of the mean approximate solution much more than costs of possible optimal solutions may differ from the cost of the mean optimal solution.

\smallskip

Also, the example we have considered shows that even though the mean optimal solution may be representative for some values of $b$, the mean approximate solution and any other possible approximate solution may be unrepresentative for much lower values of $b$. Thus, for the considered istance, while the mean optimal solution is representative for $b\!\le\!0.667$, the mean approximate solution $\tilde{x}_{\mu}$ is unrepresentative for $b\!=\!0.155$ already while the value of objective function $f(\tilde{x}_{\mu},c_{\mu})$ is closer to the maximum possible cost of a possible approximate solution rather than to its expected value.\\ 

\section{Conclusions}

The results we present show that the mean approximate solution may be unrepresentative for the set of all possible approximate solutions of IDO problem. The same is true for any other possible approximate solution: all of them may be unrepresentative. Trying to solve IDO problem and taking this into account, we must justify the representativeness of possible (optimal or approximate) solution that we obtain solving DO-problem. If we cannot give such a justification then we must solve the problem considering uncertainties in~its paremeters (e.g. costs). For example, if we need to take a solution in the situation of interval uncertainty we may obtain and analyze the problem's united solution set or its united approximate solution set and the values of cost function for them trying to predict and minimize costs.

\bigskip

The author is grateful to Omsk State University's graduate V.\,A.~Motovilov for his help that shortened the duration of computations more than $5$ times. This required near $500$ of hours of computations on PC.

\newpage

\section*{Appendix}

\subsection*{Vectors $d(m,\delta)$}

Here we present the average values for the samples
$$\{\bigl (\mathcal{P}^{(i)}(m,\delta), d^{(i)}(m,\delta)\bigr )\}_{i=1}^{100}\eqno(4)$$
for all values of $m\!=\!\overline{5,20}$, $\delta\!=\!\overline{1,5}$. 

For given $m$ and $\delta$, (4) is a sample of the pairs $\bigl (\mathcal{P}^{(i)}(m,\delta), d^{(i)}(m,\delta)\bigr )$, where
\begin{itemize}
\item[1)] $\mathcal{P}^{(i)}(m,\delta)$ is a set that contains $1000$ of ISCP instances which we generate using the presented algorithm; 
\item[2)] the vector $d^{(i)}(m,\delta)$ characterize distribution of the instances according to~probabilities of their mean approximate solutions: $$d^{(i)}(m,\delta)=\bigl (d_1^{(i)}(m,\delta),\ldots,d^{(i)}_{10}(m,\delta)\bigr ),$$
\end{itemize}
\noindent where $d^{(i)}_k(m,\delta)$ is a share of problems from $\mathcal{P}^{(i)}(m,\delta)$ (in percent) for which the probability\hspace{3pt} of\hspace{3pt} its\hspace{3pt} mean\hspace{3pt} approximate\hspace{3pt} solution\hspace{3pt} belongs to the interval $[(k\!-\!1)/10,k/10],\ k\!=\!\overline{1,10}$. 

Components of the vectors $d(m,\delta)$, $\delta\!=\!\overline{1,5}$, are placed in the rows of tables with odd numbers. Its standard deviation values are presented in accompanied tables with~even numbers.\\ 

\begin{table}[h!]
\caption{Average values of $d^{(i)}(5,\delta)$, $i\!=\!\overline{1,100}$.}
\begin{center}
\begin{tabular}{|c||c|c|c|c|c|c|c|c|c|c|}
\hline
$\delta/ b $ & $0.1$ & $0.2$ & $0.3$ & $0.4$ & $0.5$ & $0.6$ & $0.7$ & $0.8$ & $0.9$ & $1$\\
\hline
1 & 0 & 0.05 & 0.18 & 0.53 & 1.59 & 5.04 & 0.10 & 1.49 & 10.13 & 80.89 \\
\hline
2 & 0 & 0.06 & 0.35 & 1.20 & 3.64 & 4.04 & 3.65 & 7.70 & 7.39 & 71.96 \\
\hline
3 & 0 & 0.13 & 0.65 & 1.98 & 5.27 & 4.57 & 7.70 & 7.13 & 5.32 & 67.26 \\
\hline
4 & 0 & 0.20 & 1.11 & 3.31 & 6.20 & 7.39 & 7.21 & 6.17 & 5.62 & 62.77 \\
\hline
5 & 0 & 0.35 & 1.71 & 4.23 & 7.58 & 9.29 & 7.00 & 5.09 & 4.09 & 60.65 \\
\hline
\end{tabular}
\end{center}

\caption{The standard deviation of $d(5,\delta)$.}
\begin{center}
\begin{tabular}{|c||c|c|c|c|c|c|c|c|c|c|}
\hline
$\delta/ b $ & $0.1$ & $0.2$ & $0.3$ & $0.4$ & $0.5$ & $0.6$ & $0.7$ & $0.8$ & $0.9$ & $1$\\
\hline
1 & 0 & 0.07 & 0.17 & 0.23 & 0.63 & 1.17 & 0.10 & 0.45 & 2.12 & 3.94 \\
\hline
2 & 0.01 & 0.08 & 0.20 & 0.42 & 0.88 & 0.97 & 1.07 & 1.58 & 1.66 & 5.50 \\
\hline
3 & 0.01 & 0.12 & 0.27 & 0.72 & 0.97 & 1.15 & 1.53 & 1.09 & 1.10 & 5.55 \\
\hline
4 & 0.02 & 0.13 & 0.39 & 0.94 & 1.27 & 1.28 & 1.44 & 1.27 & 1.00 & 6.03 \\
\hline
5 & 0.01 & 0.18 & 0.47 & 1.13 & 1.57 & 1.51 & 1.10 & 0.95 & 0.83 & 5.92 \\
\hline
\end{tabular}
\end{center}
\end{table}

\begin{table}[h!]
\caption{Average values of $d^{(i)}(6,\delta)$, $i\!=\!\overline{1,100}$.}
\begin{center}
\begin{tabular}{|c||c|c|c|c|c|c|c|c|c|c|}
\hline
$\delta/ b $ & $0.1$ & $0.2$ & $0.3$ & $0.4$ & $0.5$ & $0.6$ & $0.7$ & $0.8$ & $0.9$ & $1$\\
\hline
1 & 0.01 & 0.06 & 0.28 & 0.61 & 2.04 & 5.83 & 0.19 & 1.83 & 11.48 & 77.67 \\
\hline
2 & 0 & 0.08 & 0.50 & 1.57 & 4.26 & 4.76 & 4.30 & 8.52 & 8.09 & 67.92 \\
\hline
3 & 0 & 0.15 & 0.95 & 2.77 & 6.23 & 5.66 & 8.83 & 7.62 & 5.83 & 61.96 \\
\hline
4 & 0 & 0.27 & 1.62 & 4.49 & 7.53 & 8.33 & 7.71 & 6.76 & 5.70 & 57.60 \\
\hline
5 & 0.01 & 0.48 & 2.65 & 5.96 & 9.42 & 10.30 & 7.60 & 5.35 & 4.31 & 53.90 \\
\hline
\end{tabular}
\end{center}

\begin{center}
\caption{The standard deviation of $d(6,\delta)$.}
\begin{tabular}{|c||c|c|c|c|c|c|c|c|c|c|}
\hline
$\delta/ b $ & $0.1$ & $0.2$ & $0.3$ & $0.4$ & $0.5$ & $0.6$ & $0.7$ & $0.8$ & $0.9$ & $1$\\
\hline
1 & 0.02 & 0.08 & 0.17 & 0.24 & 0.50 & 1.08 & 0.13 & 0.48 & 1.45 & 2.63 \\
\hline
2 & 0.02 & 0.09 & 0.22 & 0.44 & 0.82 & 0.76 & 0.84 & 1.35 & 1.40 & 4.14 \\
\hline
3 & 0.02 & 0.14 & 0.32 & 0.62 & 0.90 & 1.05 & 1.40 & 1.12 & 1.17 & 4.83 \\
\hline
4 & 0.02 & 0.18 & 0.45 & 0.80 & 1.15 & 1.12 & 1.17 & 1.05 & 0.96 & 4.66 \\
\hline
5 & 0.03 & 0.29 & 0.63 & 0.98 & 1.13 & 1.31 & 1.18 & 0.93 & 0.98 & 5.31 \\
\hline
\end{tabular}
\end{center}
\end{table}

\newpage

\begin{table}[h!]
\caption{Average values of $d^{(i)}(7,\delta)$, $i\!=\!\overline{1,100}$.}
\begin{center}
\begin{tabular}{|c||c|c|c|c|c|c|c|c|c|c|}
\hline
$\delta/ b $ & $0.1$ & $0.2$ & $0.3$ & $0.4$ & $0.5$ & $0.6$ & $0.7$ & $0.8$ & $0.9$ & $1$\\
\hline
1 & 0 & 0.04 & 0.25 & 0.57 & 2.06 & 6.02 & 0.12 & 1.85 & 11.92 & 77.16 \\
\hline
2 & 0 & 0.07 & 0.50 & 1.48 & 4.24 & 4.91 & 4.15 & 8.94 & 8.76 & 66.96 \\
\hline
3 & 0 & 0.13 & 0.96 & 2.87 & 6.30 & 5.77 & 9.66 & 8.14 & 6.08 & 60.08 \\
\hline
4 & 0.01 & 0.25 & 1.55 & 4.47 & 7.76 & 8.74 & 8.31 & 6.89 & 6.33 & 55.69 \\
\hline
5 & 0.01 & 0.45 & 2.63 & 5.74 & 9.54 & 10.55 & 8.09 & 5.87 & 5.01 & 52.07 \\
\hline
\end{tabular}
\end{center}

\caption{The standard deviation of $d(7,\delta)$.}
\begin{center}
\begin{tabular}{|c||c|c|c|c|c|c|c|c|c|c|}
\hline
$\delta/ b $ & $0.1$ & $0.2$ & $0.3$ & $0.4$ & $0.5$ & $0.6$ & $0.7$ & $0.8$ & $0.9$ & $1$\\
\hline
1 & 0.01 & 0.07 & 0.15 & 0.26 & 0.49 & 0.86 & 0.11 & 0.43 & 1.51 & 2.52 \\
\hline
2 & 0.01 & 0.09 & 0.23 & 0.38 & 0.62 & 0.92 & 0.70 & 1.22 & 1.10 & 3.54 \\
\hline
3 & 0.01 & 0.12 & 0.33 & 0.61 & 0.94 & 0.89 & 1.62 & 1.32 & 1.14 & 4.94 \\
\hline
4 & 0.04 & 0.17 & 0.37 & 0.73 & 1.13 & 1.18 & 1.43 & 1.13 & 1.25 & 5.25 \\
\hline
5 & 0.03 & 0.26 & 0.63 & 1.06 & 1.25 & 1.49 & 1.09 & 1.22 & 1.33 & 6.41 \\
\hline
\end{tabular}
\end{center}
\end{table}

\begin{table}[h!]
\caption{Average values of $d^{(i)}(8,\delta)$, $i\!=\!\overline{1,100}$.}
\begin{center}
\begin{tabular}{|c||c|c|c|c|c|c|c|c|c|c|}
\hline
$\delta/ b $ & $0.1$ & $0.2$ & $0.3$ & $0.4$ & $0.5$ & $0.6$ & $0.7$ & $0.8$ & $0.9$ & $1$\\
\hline
1 & 0 & 0.07 & 0.27 & 0.56 & 2.00 & 6.10 & 0.16 & 1.76 & 11.99 & 77.08 \\
\hline
2 & 0 & 0.09 & 0.48 & 1.47 & 3.94 & 4.91 & 4.21 & 8.90 & 8.47 & 67.53 \\
\hline
3 & 0.01 & 0.16 & 0.90 & 2.77 & 6.03 & 5.90 & 9.74 & 8.35 & 6.08 & 60.06 \\
\hline
4 & 0.02 & 0.31 & 1.56 & 4.29 & 7.26 & 8.56 & 8.21 & 7.07 & 6.76 & 55.96 \\
\hline
5 & 0.04 & 0.51 & 2.40 & 5.29 & 8.85 & 10.46 & 7.71 & 5.71 & 5.39 & 53.63 \\
\hline
\end{tabular}
\end{center}

\caption{The standard deviation of $d(8,\delta)$.}
\begin{center}
\begin{tabular}{|c||c|c|c|c|c|c|c|c|c|c|}
\hline
$\delta/ b $ & $0.1$ & $0.2$ & $0.3$ & $0.4$ & $0.5$ & $0.6$ & $0.7$ & $0.8$ & $0.9$ & $1$\\
\hline
1 & 0.01 & 0.08 & 0.16 & 0.25 & 0.45 & 0.97 & 0.12 & 0.43 & 1.46 & 2.23 \\
\hline
2 & 0 & 0.09 & 0.21 & 0.37 & 0.63 & 0.98 & 0.72 & 1.23 & 1.10 & 3.10 \\
\hline
3 & 0.03 & 0.13 & 0.30 & 0.57 & 0.75 & 0.89 & 1.33 & 1.24 & 1.17 & 4.07 \\
\hline
4 & 0.04 & 0.18 & 0.34 & 0.70 & 0.85 & 0.93 & 1.11 & 1.03 & 1.28 & 3.62 \\
\hline
5 & 0.07 & 0.22 & 0.48 & 0.74 & 1.02 & 0.87 & 0.82 & 0.74 & 0.84 & 2.01 \\
\hline
\end{tabular}
\end{center}
\end{table}

\newpage

\begin{table}[h!]
\caption{Average values of $d^{(i)}(9,\delta)$, $i\!=\!\overline{1,100}$.}
\begin{center}
\begin{tabular}{|c||c|c|c|c|c|c|c|c|c|c|}
\hline
$\delta/ b $ & $0.1$ & $0.2$ & $0.3$ & $0.4$ & $0.5$ & $0.6$ & $0.7$ & $0.8$ & $0.9$ & $1$\\
\hline
1 & 0 & 0.05 & 0.20 & 0.56 & 1.88 & 5.74 & 0.13 & 1.81 & 11.73 & 77.90 \\
\hline
2 & 0 & 0.09 & 0.41 & 1.33 & 3.91 & 4.75 & 3.81 & 8.95 & 8.67 & 68.09 \\
\hline
3 & 0 & 0.12 & 0.75 & 2.57 & 5.77 & 5.60 & 9.84 & 8.59 & 6.64 & 60.11 \\
\hline
4 & 0.01 & 0.23 & 1.34 & 3.89 & 7.10 & 8.62 & 8.84 & 7.82 & 7.48 & 54.67 \\
\hline
5 & 0.03 & 0.43 & 2.26 & 5.19 & 9.27 & 11.24 & 8.91 & 6.83 & 6.77 & 49.05 \\
\hline
\end{tabular}
\end{center}

\caption{The standard deviation of $d(9,\delta)$.}
\begin{center}
\begin{tabular}{|c||c|c|c|c|c|c|c|c|c|c|}
\hline
$\delta/ b $ & $0.1$ & $0.2$ & $0.3$ & $0.4$ & $0.5$ & $0.6$ & $0.7$ & $0.8$ & $0.9$ & $1$\\
\hline
1 & 0.01 & 0.08 & 0.14 & 0.24 & 0.41 & 0.75 & 0.11 & 0.43 & 0.91 & 1.31 \\
\hline
2 & 0.02 & 0.10 & 0.19 & 0.38 & 0.73 & 0.76 & 0.64 & 0.98 & 0.93 & 1.42 \\
\hline
3 & 0.02 & 0.12 & 0.29 & 0.60 & 0.77 & 0.78 & 1.00 & 1.09 & 1.06 & 2.46 \\
\hline
4 & 0.03 & 0.16 & 0.39 & 0.62 & 0.79 & 0.88 & 0.89 & 1.10 & 1.49 & 2.35 \\
\hline
5 & 0.05 & 0.21 & 0.46 & 0.73 & 0.96 & 1.39 & 1.12 & 1.28 & 1.38 & 4.13 \\
\hline
\end{tabular}
\end{center}
\end{table}

\begin{table}[h!]
\caption{Average values of $d^{(i)}(10,\delta)$, $i\!=\!\overline{1,100}$.}
\begin{center}
\begin{tabular}{|c||c|c|c|c|c|c|c|c|c|c|}
\hline
$\delta/ b $ & $0.1$ & $0.2$ & $0.3$ & $0.4$ & $0.5$ & $0.6$ & $0.7$ & $0.8$ & $0.9$ & $1$\\
\hline
1 & 0 & 0.05 & 0.23 & 0.48 & 1.72 & 6.08 & 0.13 & 1.82 & 11.87 & 77.61 \\
\hline
2 & 0 & 0.07 & 0.42 & 1.25 & 3.55 & 5.09 & 3.73 & 9.48 & 8.84 & 67.57 \\
\hline
3 & 0.01 & 0.15 & 0.81 & 2.31 & 5.47 & 5.67 & 9.90 & 8.82 & 6.92 & 59.94 \\
\hline
4 & 0.01 & 0.23 & 1.32 & 3.84 & 7.09 & 8.57 & 9.09 & 8.46 & 8.87 & 52.53 \\
\hline
5 & 0.02 & 0.40 & 2.10 & 4.95 & 9.17 & 11.36 & 9.21 & 7.84 & 8.12 & 46.81 \\
\hline
\end{tabular}
\end{center}

\caption{Standard deviation of $d(10,\delta)$.}
\begin{center}
\begin{tabular}{|c||c|c|c|c|c|c|c|c|c|c|}
\hline
$\delta/ b $ & $0.1$ & $0.2$ & $0.3$ & $0.4$ & $0.5$ & $0.6$ & $0.7$ & $0.8$ & $0.9$ & $1$\\
\hline
1 & 0.01 & 0.07 & 0.15 & 0.21 & 0.45 & 0.88 & 0.12 & 0.42 & 1.11 & 1.46 \\
\hline
2 & 0 & 0.08 & 0.20 & 0.36 & 0.69 & 0.79 & 0.71 & 1.00 & 1.01 & 1.72 \\
\hline
3 & 0.03 & 0.15 & 0.31 & 0.54 & 0.80 & 0.77 & 0.98 & 0.82 & 0.94 & 1.65 \\
\hline
4 & 0.03 & 0.17 & 0.42 & 0.69 & 0.83 & 0.91 & 1.08 & 1.34 & 1.20 & 2.65 \\
\hline
5 & 0.04 & 0.20 & 0.42 & 0.79 & 0.88 & 1.18 & 1.17 & 1.07 & 1.53 & 3.54 \\
\hline
\end{tabular}
\end{center}
\end{table}

\newpage

\begin{table}[h!]
\caption{Average values of $d^{(i)}(11,\delta)$, $i\!=\!\overline{1,100}$.}
\begin{center}
\begin{tabular}{|c||c|c|c|c|c|c|c|c|c|c|}
\hline
$\delta/ b $ & $0.1$ & $0.2$ & $0.3$ & $0.4$ & $0.5$ & $0.6$ & $0.7$ & $0.8$ & $0.9$ & $1$\\
\hline
1 & 0 & 0.04 & 0.18 & 0.43 & 1.53 & 5.71 & 0.16 & 1.98 & 11.55 & 78.43 \\
\hline
2 & 0 & 0.06 & 0.41 & 1.24 & 3.49 & 5.09 & 3.99 & 9.16 & 8.84 & 67.72 \\
\hline
3 & 0 & 0.13 & 0.81 & 2.26 & 5.48 & 5.87 & 9.73 & 8.88 & 7.46 & 59.38 \\
\hline
4 & 0.01 & 0.23 & 1.19 & 3.75 & 6.88 & 8.91 & 9.29 & 8.47 & 9.56 & 51.71 \\
\hline
5 & 0.02 & 0.40 & 2.13 & 4.94 & 9.09 & 11.38 & 9.38 & 8.00 & 8.15 & 46.49 \\
\hline
\end{tabular}
\end{center}

\caption{The standard deviation of $d(11,\delta)$.}
\begin{center}
\begin{tabular}{|c||c|c|c|c|c|c|c|c|c|c|}
\hline
$\delta/ b $ & $0.1$ & $0.2$ & $0.3$ & $0.4$ & $0.5$ & $0.6$ & $0.7$ & $0.8$ & $0.9$ & $1$\\
\hline
1 & 0 & 0.07 & 0.13 & 0.23 & 0.46 & 0.83 & 0.14 & 0.47 & 1.01 & 1.34 \\
\hline
2 & 0.02 & 0.08 & 0.19 & 0.36 & 0.71 & 0.65 & 0.68 & 0.87 & 1.02 & 1.40 \\
\hline
3 & 0.02 & 0.13 & 0.30 & 0.51 & 0.85 & 0.80 & 0.99 & 0.94 & 0.97 & 1.76 \\
\hline
4 & 0.04 & 0.14 & 0.37 & 0.71 & 0.94 & 0.88 & 1.05 & 1.19 & 1.54 & 2.28 \\
\hline
5 & 0.05 & 0.19 & 0.49 & 0.74 & 0.78 & 1.27 & 1.08 & 1.01 & 1.63 & 2.97 \\
\hline
\end{tabular}
\end{center}
\end{table}

\begin{table}[h!]
\caption{Average values of $d^{(i)}(13,\delta)$, $i\!=\!\overline{1,100}$.}
\begin{center}
\begin{tabular}{|c||c|c|c|c|c|c|c|c|c|c|}
\hline
$\delta/ b $ & $0.1$ & $0.2$ & $0.3$ & $0.4$ & $0.5$ & $0.6$ & $0.7$ & $0.8$ & $0.9$ & $1$\\
\hline
1 & 0 & 0.05 & 0.20 & 0.50 & 1.71 & 5.82 & 0.54 & 1.92 & 11.27 & 77.98 \\
\hline
2 & 0 & 0.08 & 0.39 & 1.27 & 3.48 & 5.19 & 4.14 & 9.32 & 9.16 & 66.97 \\
\hline
3 & 0 & 0.14 & 0.76 & 2.32 & 5.45 & 6.20 & 10.03 & 9.01 & 7.78 & 58.32 \\
\hline
4 & 0.01 & 0.27 & 1.35 & 3.86 & 7.05 & 8.73 & 9.23 & 8.72 & 9.15 & 51.63 \\
\hline
5 & 0.01 & 0.42 & 2.18 & 4.83 & 9.03 & 11.58 & 9.62 & 8.56 & 8.71 & 45.01 \\
\hline
\end{tabular}
\end{center}

\caption{The standard deviation of $d(12,\delta)$.}
\begin{center}
\begin{tabular}{|c||c|c|c|c|c|c|c|c|c|c|}
\hline
$\delta/ b $ & $0.1$ & $0.2$ & $0.3$ & $0.4$ & $0.5$ & $0.6$ & $0.7$ & $0.8$ & $0.9$ & $1$\\
\hline
1 & 0 & 0.08 & 0.14 & 0.25 & 0.48 & 0.86 & 0.26 & 0.49 & 0.95 & 1.10 \\
\hline
2 & 0.01 & 0.10 & 0.20 & 0.39 & 0.67 & 0.7 & 0.71 & 0.79 & 0.93 & 1.55 \\
\hline
3 & 0.01 & 0.12 & 0.32 & 0.59 & 0.95 & 0.77 & 0.97 & 0.92 & 0.93 & 1.50 \\
\hline
4 & 0.02 & 0.17 & 0.39 & 0.66 & 0.84 & 0.94 & 1.02 & 1.01 & 1.37 & 2.30 \\
\hline
5 & 0.04 & 0.25 & 0.56 & 0.91 & 0.94 & 1.19 & 1.05 & 1.10 & 1.60 & 2.76 \\
\hline
\end{tabular}
\end{center}
\end{table}

\newpage

\begin{table}[h!]
\caption{Average values of $d^{(i)}(13,\delta)$, $i\!=\!\overline{1,100}$.}
\begin{center}
\begin{tabular}{|c||c|c|c|c|c|c|c|c|c|c|}
\hline
$\delta/ b $ & $0.1$ & $0.2$ & $0.3$ & $0.4$ & $0.5$ & $0.6$ & $0.7$ & $0.8$ & $0.9$ & $1$\\
\hline
1 & 0 & 0.04 & 0.22 & 0.41 & 1.63 & 6.02 & 0.42 & 1.63 & 11.56 & 78.06 \\
\hline
2 & 0 & 0.05 & 0.38 & 1.09 & 3.38 & 5.35 & 3.87 & 9.84 & 9.63 & 66.41 \\
\hline
3 & 0 & 0.12 & 0.68 & 2.20 & 5.29 & 6.20 & 10.3 & 9.43 & 7.96 & 57.82 \\
\hline
4 & 0 & 0.19 & 1.16 & 3.35 & 6.79 & 8.32 & 9.68 & 9.24 & 9.82 & 51.45 \\
\hline
5 & 0.01 & 0.36 & 1.87 & 4.45 & 8.6 & 11.86 & 9.79 & 8.27 & 9.58 & 45.18 \\
\hline
\end{tabular}
\end{center}

\caption{The standard deviation of $d(13,\delta)$.}
\begin{center}
\begin{tabular}{|c||c|c|c|c|c|c|c|c|c|c|}
\hline
$\delta/ b $ & $0.1$ & $0.2$ & $0.3$ & $0.4$ & $0.5$ & $0.6$ & $0.7$ & $0.8$ & $0.9$ & $1$\\
\hline
1 & 0.02 & 0.06 & 0.13 & 0.21 & 0.44 & 0.79 & 0.2 & 0.44 & 1.05 & 1.65 \\
\hline
2 & 0.01 & 0.08 & 0.21 & 0.37 & 0.74 & 0.72 & 0.68 & 0.92 & 0.92 & 2.08 \\
\hline
3 & 0.01 & 0.11 & 0.31 & 0.59 & 0.92 & 0.79 & 0.89 & 0.98 & 0.95 & 2.17 \\
\hline
4 & 0.02 & 0.17 & 0.46 & 0.91 & 1.23 & 1.06 & 0.99 & 1.09 & 1.66 & 3.22 \\
\hline
5 & 0.03 & 0.25 & 0.59 & 1.21 & 1.61 & 1.08 & 1.14 & 1.14 & 1.74 & 4.59 \\
\hline
\end{tabular}
\end{center}
\end{table}

\begin{table}[h!]
\caption{Average values of $d^{(i)}(14,\delta)$, $i\!=\!\overline{1,100}$.}
\begin{center}
\begin{tabular}{|c||c|c|c|c|c|c|c|c|c|c|}
\hline
$\delta/ b $ & $0.1$ & $0.2$ & $0.3$ & $0.4$ & $0.5$ & $0.6$ & $0.7$ & $0.8$ & $0.9$ & $1$\\
\hline
1 & 0 & 0.03 & 0.19 & 0.38 & 1.42 & 5.11 & 0.87 & 2.02 & 10.41 & 79.57 \\
\hline
2 & 0 & 0.04 & 0.32 & 1.13 & 2.98 & 5.08 & 4.04 & 9.07 & 8.96 & 68.38 \\
\hline
3 & 0 & 0.10 & 0.69 & 2.22 & 4.99 & 6.49 & 9.72 & 8.72 & 8.31 & 58.75 \\
\hline
4 & 0.01 & 0.23 & 1.25 & 3.60 & 6.95 & 8.81 & 9.71 & 9.04 & 10.21 & 50.17 \\
\hline
5 & 0.02 & 0.39 & 2.14 & 4.72 & 8.80 & 11.73 & 10 & 8.96 & 9.94 & 43.29 \\
\hline
\end{tabular}
\end{center}

\caption{The standard deviation of $d(14,\delta)$.}
\begin{center}
\begin{tabular}{|c||c|c|c|c|c|c|c|c|c|c|}
\hline
$\delta/ b $ & $0.1$ & $0.2$ & $0.3$ & $0.4$ & $0.5$ & $0.6$ & $0.7$ & $0.8$ & $0.9$ & $1$\\
\hline
1 & 0 & 0.05 & 0.14 & 0.23 & 0.52 & 0.87 & 0.31 & 0.53 & 1.38 & 2.52 \\
\hline
2 & 0.01 & 0.07 & 0.20 & 0.47 & 0.94 & 0.69 & 1.01 & 1.13 & 0.95 & 3.53 \\
\hline
3 & 0.02 & 0.10 & 0.36 & 0.68 & 1.11 & 0.89 & 1.29 & 1.05 & 0.97 & 3.97 \\
\hline
4 & 0.03 & 0.15 & 0.50 & 0.94 & 1.35 & 1.23 & 1.05 & 1.11 & 1.12 & 4.36 \\
\hline
5 & 0.05 & 0.23 & 0.69 & 1.21 & 1.73 & 1.17 & 1.13 & 1.11 & 1.04 & 4.86 \\
\hline
\end{tabular}
\end{center}
\end{table}

\newpage

\begin{table}[h!]
\caption{Average values of $d^{(i)}(15,\delta)$, $i\!=\!\overline{1,100}$.}
\begin{center}
\begin{tabular}{|c||c|c|c|c|c|c|c|c|c|c|}
\hline
$\delta/ b $ & $0.1$ & $0.2$ & $0.3$ & $0.4$ & $0.5$ & $0.6$ & $0.7$ & $0.8$ & $0.9$ & $1$\\
\hline
1 & 0 & 0.03 & 0.22 & 0.39 & 1.39 & 5.85 & 0.71 & 1.70 & 11.54 & 78.16 \\
\hline
2 & 0 & 0.05 & 0.35 & 1.08 & 3.07 & 5.46 & 4.06 & 9.28 & 9.45 & 67.19 \\
\hline
3 & 0 & 0.10 & 0.60 & 1.88 & 4.83 & 6.10 & 10.31 & 9.42 & 8.47 & 58.29 \\
\hline
4 & 0.01 & 0.18 & 1.08 & 3.32 & 6.91 & 8.72 & 9.91 & 9.64 & 10.54 & 49.67 \\
\hline
5 & 0.01 & 0.35 & 1.84 & 4.53 & 8.73 & 12.26 & 10.18 & 9.11 & 10.96 & 41.99 \\
\hline
\end{tabular}
\end{center}

\caption{The standard deviation of $d(15,\delta)$.}
\begin{center}
\begin{tabular}{|c||c|c|c|c|c|c|c|c|c|c|}
\hline
$\delta/ b $ & $0.1$ & $0.2$ & $0.3$ & $0.4$ & $0.5$ & $0.6$ & $0.7$ & $0.8$ & $0.9$ & $1$\\
\hline
1 & 0 & 0.05 & 0.14 & 0.22 & 0.47 & 0.92 & 0.22 & 0.46 & 1.39 & 2.55 \\
\hline
2 & 0.01 & 0.09 & 0.25 & 0.41 & 0.85 & 0.87 & 0.99 & 1.11 & 1.15 & 3.98 \\
\hline
3 & 0.01 & 0.10 & 0.33 & 0.79 & 1.14 & 0.86 & 1.18 & 1.36 & 0.87 & 4.16 \\
\hline
4 & 0.03 & 0.15 & 0.51 & 1.07 & 1.32 & 1.26 & 1.11 & 1.06 & 1.73 & 3.59 \\
\hline
5 & 0.03 & 0.26 & 0.65 & 1.20 & 1.68 & 1.00 & 1.20 & 1.19 & 1.70 & 4.31 \\
\hline
\end{tabular}
\end{center}
\end{table}

\begin{table}[h!]
\caption{Average values of $d^{(i)}(16,\delta)$, $i\!=\!\overline{1,100}$.}
\begin{center}
\begin{tabular}{|c||c|c|c|c|c|c|c|c|c|c|}
\hline
$\delta/ b $ & $0.1$ & $0.2$ & $0.3$ & $0.4$ & $0.5$ & $0.6$ & $0.7$ & $0.8$ & $0.9$ & $1$\\
\hline
1 & 0 & 0.03 & 0.15 & 0.32 & 1.23 & 5.20 & 1.26 & 1.97 & 9.90 & 79.93 \\
\hline
2 & 0 & 0.05 & 0.28 & 0.98 & 2.72 & 5.55 & 4.25 & 9.15 & 9.17 & 67.85 \\
\hline
3 & 0 & 0.08 & 0.48 & 1.71 & 4.41 & 6.83 & 9.62 & 9.38 & 9.17 & 58.31 \\
\hline
4 & 0 & 0.13 & 0.90 & 2.89 & 6.24 & 8.65 & 10.00 & 9.77 & 11.56 & 49.86 \\
\hline
5 & 0.01 & 0.28 & 1.60 & 4.14 & 8.12 & 12.06 & 10.37 & 10.21 & 11.94 & 41.26 \\
\hline
\end{tabular}
\end{center}

\caption{The standard deviation of $d(16,\delta)$.}
\begin{center}
\begin{tabular}{|c||c|c|c|c|c|c|c|c|c|c|}
\hline
$\delta/ b $ & $0.1$ & $0.2$ & $0.3$ & $0.4$ & $0.5$ & $0.6$ & $0.7$ & $0.8$ & $0.9$ & $1$\\
\hline
1 & 0 & 0.05 & 0.14 & 0.24 & 0.51 & 0.79 & 0.45 & 0.54 & 1.49 & 2.81 \\
\hline
2 & 0 & 0.07 & 0.19 & 0.40 & 0.79 & 0.88 & 0.88 & 1.14 & 0.98 & 3.06 \\
\hline
3 & 0.02 & 0.10 & 0.32 & 0.68 & 1.26 & 0.84 & 0.95 & 0.98 & 1.22 & 2.87 \\
\hline
4 & 0.01 & 0.12 & 0.52 & 1.21 & 1.65 & 1.38 & 1.15 & 1.16 & 1.73 & 4.71 \\
\hline
5 & 0.03 & 0.24 & 0.67 & 1.26 & 1.67 & 1.21 & 1.15 & 1.22 & 1.36 & 4.94 \\
\hline
\end{tabular}
\end{center}
\end{table}

\newpage

\begin{table}[h!]
\caption{Average values of $d^{(i)}(17,\delta)$, $i\!=\!\overline{1,100}$.}
\begin{center}
\begin{tabular}{|c||c|c|c|c|c|c|c|c|c|c|}
\hline
$\delta/ b $ & $0.1$ & $0.2$ & $0.3$ & $0.4$ & $0.5$ & $0.6$ & $0.7$ & $0.8$ & $0.9$ & $1$\\
\hline
1 & 0 & 0.01 & 0.17 & 0.31 & 1.17 & 5.71 & 0.90 & 1.65 & 10.49 & 79.58 \\
\hline
2 & 0 & 0.03 & 0.26 & 0.97 & 2.77 & 5.59 & 3.79 & 9.69 & 9.58 & 67.31 \\
\hline
3 & 0 & 0.08 & 0.56 & 2.01 & 5.06 & 6.59 & 10.51 & 9.69 & 8.78 & 56.71 \\
\hline
4 & 0 & 0.14 & 1.02 & 3.20 & 6.83 & 8.85 & 10.26 & 10.08 & 11.64 & 47.97 \\
\hline
5 & 0.01 & 0.30 & 1.72 & 4.43 & 8.78 & 12.52 & 10.82 & 10.00 & 11.83 & 39.59 \\
\hline
\end{tabular}
\end{center}

\caption{The standard deviation of $d(17,\delta)$.}
\begin{center}
\begin{tabular}{|c||c|c|c|c|c|c|c|c|c|c|}
\hline
$\delta/ b $ & $0.1$ & $0.2$ & $0.3$ & $0.4$ & $0.5$ & $0.6$ & $0.7$ & $0.8$ & $0.9$ & $1$\\
\hline
1 & 0.01 & 0.04 & 0.15 & 0.20 & 0.54 & 1.02 & 0.32 & 0.49 & 1.58 & 3.00 \\
\hline
2 & 0 & 0.05 & 0.18 & 0.42 & 0.86 & 0.77 & 0.84 & 1.26 & 1.05 & 3.76 \\
\hline
3 & 0.01 & 0.09 & 0.29 & 0.67 & 1.20 & 1.07 & 1.44 & 1.01 & 1.17 & 4.38 \\
\hline
4 & 0.02 & 0.13 & 0.48 & 1.01 & 1.47 & 1.38 & 1.11 & 1.25 & 1.20 & 4.82 \\
\hline
5 & 0.02 & 0.22 & 0.69 & 1.37 & 1.85 & 1.09 & 1.08 & 1.09 & 1.33 & 5.06 \\
\hline
\end{tabular}
\end{center}
\end{table}

\begin{table}[h!]
\caption{Average values of $d^{(i)}(18,\delta)$, $i\!=\!\overline{1,100}$.}
\begin{center}
\begin{tabular}{|c||c|c|c|c|c|c|c|c|c|c|}
\hline
$\delta/ b $ & $0.1$ & $0.2$ & $0.3$ & $0.4$ & $0.5$ & $0.6$ & $0.7$ & $0.8$ & $0.9$ & $1$\\
\hline
1 & 0 & 0.02 & 0.17 & 0.35 & 1.20 & 5.07 & 1.43 & 2.09 & 10.30 & 79.35 \\
\hline
2 & 0 & 0.03 & 0.31 & 0.93 & 2.69 & 5.48 & 4.42 & 8.78 & 9.59 & 67.78 \\
\hline
3 & 0 & 0.09 & 0.53 & 1.84 & 4.80 & 6.79 & 10.04 & 9.75 & 9.46 & 56.71 \\
\hline
4 & 0.01 & 0.14 & 0.98 & 3.11 & 6.68 & 9.62 & 10.41 & 10.53 & 11.33 & 47.20 \\
\hline
5 & 0.01 & 0.27 & 1.73 & 4.39 & 8.73 & 12.60 & 11.04 & 10.71 & 12.10 & 38.40 \\
\hline
\end{tabular}
\end{center}

\caption{The standard deviation of $d(18,\delta)$.}
\begin{center}
\begin{tabular}{|c||c|c|c|c|c|c|c|c|c|c|}
\hline
$\delta/ b $ & $0.1$ & $0.2$ & $0.3$ & $0.4$ & $0.5$ & $0.6$ & $0.7$ & $0.8$ & $0.9$ & $1$\\
\hline
1 & 0 & 0.05 & 0.13 & 0.20 & 0.46 & 0.80 & 0.44 & 0.5 & 1.24 & 2.59 \\
\hline
2 & 0 & 0.05 & 0.22 & 0.42 & 0.80 & 0.82 & 0.88 & 1.11 & 1.34 & 3.69 \\
\hline
3 & 0.01 & 0.11 & 0.29 & 0.57 & 1.30 & 0.98 & 1.27 & 1.02 & 1.25 & 4.28 \\
\hline
4 & 0.02 & 0.11 & 0.52 & 0.87 & 1.36 & 1.25 & 1.17 & 1.06 & 1.15 & 4.79 \\
\hline
5 & 0.03 & 0.19 & 0.72 & 1.29 & 1.72 & 1.16 & 1.19 & 1.16 & 1.36 & 5.46 \\
\hline
\end{tabular}
\end{center}
\end{table}

\newpage

\begin{table}[h!]
\caption{Average values of $d^{(i)}(19,\delta)$, $i\!=\!\overline{1,100}$.}
\begin{center}
\begin{tabular}{|c||c|c|c|c|c|c|c|c|c|c|}
\hline
$\delta/ b $ & $0.1$ & $0.2$ & $0.3$ & $0.4$ & $0.5$ & $0.6$ & $0.7$ & $0.8$ & $0.9$ & $1$\\
\hline
1 & 0 & 0.02 & 0.13 & 0.32 & 1.21 & 4.95 & 1.78 & 2.58 & 9.52 & 79.46 \\
\hline
2 & 0 & 0.03 & 0.26 & 0.84 & 2.57 & 5.58 & 4.46 & 9.18 & 9.14 & 67.93 \\
\hline
3 & 0 & 0.05 & 0.44 & 1.65 & 4.46 & 7.09 & 10.06 & 9.35 & 10.18 & 56.72 \\
\hline
4 & 0 & 0.11 & 0.88 & 2.79 & 6.53 & 9.55 & 10.72 & 10.53 & 12.54 & 46.33 \\
\hline
5 & 0 & 0.22 & 1.45 & 3.89 & 8.42 & 12.75 & 11.60 & 11.83 & 12.68 & 37.15 \\
\hline
\end{tabular}
\end{center}

\caption{The standard deviation of $d(19,\delta)$.}
\begin{center}
\begin{tabular}{|c||c|c|c|c|c|c|c|c|c|c|}
\hline
$\delta/ b $ & $0.1$ & $0.2$ & $0.3$ & $0.4$ & $0.5$ & $0.6$ & $0.7$ & $0.8$ & $0.9$ & $1$\\
\hline
1 & 0 & 0.05 & 0.12 & 0.20 & 0.41 & 0.77 & 0.49 & 0.66 & 1.17 & 2.59 \\
\hline
2 & 0.01 & 0.06 & 0.17 & 0.38 & 0.86 & 0.81 & 1.05 & 1.20 & 0.98 & 3.50 \\
\hline
3 & 0.01 & 0.07 & 0.24 & 0.63 & 1.29 & 1.16 & 1.23 & 1.09 & 1.00 & 4.15 \\
\hline
4 & 0.02 & 0.12 & 0.41 & 0.87 & 1.49 & 1.35 & 1.00 & 1.17 & 1.25 & 4.67 \\
\hline
5 & 0.01 & 0.16 & 0.57 & 1.14 & 1.80 & 1.37 & 1.32 & 0.94 & 1.36 & 5.07 \\
\hline
\end{tabular}
\end{center}
\end{table}

\begin{table}[h!]
\caption{Average values of $d^{(i)}(20,\delta)$, $i\!=\!\overline{1,100}$.}
\begin{center}
\begin{tabular}{|c||c|c|c|c|c|c|c|c|c|c|}
\hline
$\delta/ b $ & $0.1$ & $0.2$ & $0.3$ & $0.4$ & $0.5$ & $0.6$ & $0.7$ & $0.8$ & $0.9$ & $1$\\
\hline
1 & 0 & 0.01 & 0.10 & 0.29 & 1.01 & 5.35 & 1.42 & 1.89 & 9.92 & 80.00 \\
\hline
2 & 0 & 0.04 & 0.19 & 0.79 & 2.43 & 5.69 & 4.12 & 9.37 & 9.46 & 67.91 \\
\hline
3 & 0 & 0.06 & 0.46 & 1.67 & 4.64 & 7.04 & 10.71 & 10.02 & 10.55 & 54.84 \\
\hline
4 & 0 & 0.11 & 0.78 & 2.74 & 6.46 & 9.48 & 10.87 & 11.10 & 12.94 & 45.52 \\
\hline
5 & 0.01 & 0.20 & 1.32 & 3.98 & 8.62 & 12.92 & 11.89 & 11.46 & 12.97 & 36.63 \\
\hline
\end{tabular}
\end{center}

\caption{The standard deviation of $d(20,\delta)$.}
\begin{center}
\begin{tabular}{|c||c|c|c|c|c|c|c|c|c|c|}
\hline
$\delta/ b $ & $0.1$ & $0.2$ & $0.3$ & $0.4$ & $0.5$ & $0.6$ & $0.7$ & $0.8$ & $0.9$ & $1$\\
\hline
1 & 0.01 & 0.03 & 0.11 & 0.22 & 0.40 & 0.91 & 0.37 & 0.37 & 1.72 & 2.58 \\
\hline
2 & 0 & 0.06 & 0.15 & 0.33 & 0.83 & 0.76 & 0.71 & 1.31 & 1.08 & 3.50 \\
\hline
3 & 0 & 0.07 & 0.25 & 0.57 & 1.32 & 0.90 & 1.31 & 0.98 & 0.88 & 3.64 \\
\hline
4 & 0.01 & 0.12 & 0.42 & 0.95 & 1.65 & 1.52 & 1.29 & 1.15 & 1.13 & 5.66 \\
\hline
5 & 0.02 & 0.15 & 0.58 & 1.25 & 2.09 & 1.44 & 1.40 & 1.24 & 1.31 & 5.85 \\
\hline
\end{tabular}
\end{center}
\end{table}

\end{document}